\documentclass[12pt]{amsart}
\usepackage{amsmath,amscd,amssymb,amsfonts}
\theoremstyle{plain}
\newtheorem{theorem}{Theorem}
\newtheorem{lem}[theorem]{Lemma}
\newtheorem{question}[theorem]{Question}
\newtheorem{definition}[theorem]{Definition}
\newtheorem{example}[theorem]{Example}
\numberwithin{theorem}{section}
\numberwithin{equation}{section}

\newcommand{\Hom}{{\rm Hom}}

\newcommand{\sD}{{\mathcal D}}

\newcommand{\sL}{{\mathcal L}}

\newcommand{\sO}{{\mathcal O}}

\newcommand{\sT}{{\mathcal T}}

\newcommand{\C}{{\mathbb C}}

\renewcommand{\P}{{\mathbb P}}
\newcommand{\Q}{{\mathbb Q}}
\newcommand{\R}{{\mathbb R}}

\newcommand{\Z}{{\mathbb Z}}

\begin{document}

\title[Picard-Fuchs equations, Integrable Systems and ...]
{Picard-Fuchs equations, Integrable Systems and higher
algebraic K-theory}
 
\author{Pedro Luis del Angel}
\address{Pedro Luis del Angel, Cimat, Guanajuato, Mexico}

\author{Stefan M\"uller-Stach}
\address{Stefan M\"uller-Stach, McMaster University, Canada}

\date{}
\begin{abstract} 
This paper continues the work done in \cite{dAMS} and is
an attempt to establish a conceptual framework which generalizes the work 
of Manin \cite{M} on the relation between non-linear second order ODE of 
type Painlev\'e VI and integrable systems. The principle behind everything 
is a strong interaction between K-theory
and Picard-Fuchs type differential equations via Abel-Jacobi maps.
Our main result is an extension of a theorem of Donagi and Markman \cite{DM1}.
\end{abstract}
\thanks{We thank DFG, Schwerpunkt and Heisenberg program, 
Conacyt (grant 37557-E), CIMAT, Fields Institute, McMaster University 
and Universit\"at Essen for supporting this project}
\maketitle

\begin{center}{\it Dedicated to Andrei Tyurin}\end{center}

\section*{Introduction}

In \cite{M}, building up on his work on the functional Mordell conjecture, 
Y.I.Manin has found a framework in which the relation between 
Picard-Fuchs differential equations for the Legendre family of elliptic curves
and the non-linear equations of type
Painlev\'e VI due to Richard Fuchs can be connected to mathematical physics
and the theory of integrable systems. One consequence of this is an approach 
for the understanding of mirror symmetry in the case of Fano manifolds. \\
Inspired by his work and the related work of Griffiths \cite{G4} about 
differential equations satisfied by normal functions associated to classical 
cycles, we have imitated the relation between periods and 
non-linear equations in the case of K3 surfaces in \cite{dAMS}. It turned out 
that there the role of classical algebraic cycles, e.g. sums of points, has 
to be replaced by elements in higher algebraic K-groups \cite{Q}, resp. 
motivic cohomology groups \cite{B2}.\\
In this paper we attempt to push these ideas further in the case of families 
of Calabi-Yau manifolds $f: X \to B$ of relative dimension $d$. We study the 
differential equations obtained from algebraic cycles in higher Chow groups 
$CH^p(X,n)$ which have good intersection with each fiber of $f$ 
when $2p-n-1=d$. These can be looked at as generalizations of non-linear 
ODE of order less than or equal to 
the rank of the local system associated to $R^df_*\C$. 
We describe the connection to the theory of analytically completely integrable 
Hamiltonian systems (with non-complete fibers) and give several examples where
the computations are partially worked out. This part builds up on work of Donagi 
and Markman \cite{DM1}.
The upshot is that algebraic K-theory brings
a new flavour to the theory of Hamiltonian systems and conversely we get new 
insight into differential integrable systems via generalized Abel-Jacobi 
maps. There is a general understanding that
mirror symmetry also plays a role in such investigations, but we are still 
far from understanding the connection to mathematical physics. \\
We want to emphasize that these notes have a certain survey character 
required by the publisher. However all new results come with a self--contained 
proof, except some announcements made in the last section. 

\section{Symbol map and Picard-Fuchs operators}
\label{picardfuchs}

Let $f: X \to B$ a smooth family of projective manifolds of dimension $d$ 
over $\C$, where we assume that $B$ is a smooth curve. Recall that the short 
exact sequence
$$ 
0 \to \sT_{X/B} \to \sT_X \to f^*\sT_B\to 0 
$$
induces the Kodaira-Spencer map: 
$$
\sT_B = f_*f^*\sT_B \to R^1f_*\sT_{X/B}
$$
whose iterated $d$-fold cup product (a standard method in 
IVHS) produces a map
$$
S^d\sT_B \to S^d R^1f_*\sT_{X/B} 
\to \Hom(f_*\Omega^d_{X/B}, R^df_*\sO) \cong (\sL^\vee)^2
$$
since $ (R^df_*\sO)^\vee = f_*\Omega^d_{X/B} =: \sL$, a vector bundle over 
$B$, which is a line bundle in the case where all fibers have trivial 
canonical bundle, e.g. in the case of Calabi-Yau manifolds. Summarizing, 
we obtain a symbol map
$$
S^d\sT_B\to (\sL^\vee)^2,
$$
or, what amounts to the same, a section of $S^d\Omega^1_B\otimes(\sL^\vee)^2$. 
It is also possible to write a version of this with logarithmic poles along the
boundary divisor of a smooth compactification of B as follows:
$$
S^d\sT_B(-\log \Sigma)\to (\bar{\sL}^\vee)^2,
$$
where one has the following commutative diagram
$$
\begin{matrix}
\Delta & \subset & \bar X & \leftarrow & X \cr
\downarrow & & \hskip4mm \downarrow \bar f & & \hskip4mm  \downarrow f  \cr
\Sigma & \subset & \bar B & \leftarrow &B \end{matrix}
$$
and $\bar\sL := \bar f_*\Omega^d_{X/B}(\log\Delta)$. Here $\Sigma = \bar B 
\setminus B$ and $\Delta=\bar f^{-1}(\Sigma)$ are reduced normal crossing 
divisors. This situation can be obtained after semistable reduction.   \\
For a smooth family as before $f: X \to B$, one has a canonical 
class of local systems on $B$, 
given by the vanishing of the Gau\ss-Manin connection $\nabla$ 
on the higher direct images of $\C$,
$$
\nabla: R^kf_*\C\otimes\sO_B \to R^kf_*\C \otimes \Omega^1_B.
$$
Let $\sD$ be the sheaf of algebraic differential operators 
on $B$. Remember that $\sD$ admits a filtration by the order of the operator, 
in particular one has a short exact sequence
$$ 
0\to \sD^{\le m-1}\to \sD^{\le m}\to S^m\sT_B \to 0 
$$
A homogenous element in the image of $\sD^{\le m}$ 
under the last map will be called an $m$-{\it symbol}. Differential operators 
act via the Lie derivative on relative holomorphic $k$--forms:
$$
\sD \otimes f_* \Omega^k_{X/B} \to R^kf_* \C \otimes \sO_B, \quad \sigma 
\otimes \omega \mapsto {\mathcal L}_\sigma(\omega),
$$
see \cite{L}. If the rank of the (irreducible) local system  
$R^kf_*\C$ is $m$, to any $m$-symbol $\sigma$
one can associate a {\it unique} differential operator $D_\sigma$, called the 
{\it Picard-Fuchs operator} of the symbol $\sigma$, 
which annihilates all sections of $R^kf_*\C$ (\cite{Bo}, Corollary III.2.3.2). 
The associated homogeneous differential equation $D_\sigma \varphi=0$ 
is called the {\it Picard-Fuchs equation}. In particular every such  $D_\sigma$
is contained in the left kernel of 
$\sD \otimes f_* \Omega^k_{X/B} \to R^kf_* \C \otimes \sO_B$. The collection 
of all $D_\sigma$ defines an ideal in the sheaf $\sD$ of algebraic 
differential operators. In other words, for every holomorphic relative 
$k$--form $\omega$
and every $m$-symbol $\sigma$ there is a $(k-1)$-form $\beta$ such that
$$
D_\sigma \omega = d_{\rm rel} \beta, 
$$
which we call the  {\it inhomogenous Picard-Fuchs equation}.
More generally, if $\R^kf_*\C$ is not an irreducible representation of the 
monodromy group, let $r$ be the rank of a maximal
irreducible representation $M$ contained in it. Then to any $r$-symbol $\sigma$ 
there corresponds a unique differential operator $D_\sigma$ annihilating all
sections of $M$. Sometimes, if we do not want to specify $\sigma$ we will
write $D_{\rm PF}$ instead of $D_\sigma$. But the reader should be warned that
these differential ideals are not always principal ideals.
After integration over non-closed topological chains, the equation 
$D_\sigma \omega = d_{\rm rel} \beta$ is mainly responsible
for the appearance of new differential equations. We will use this
construction in the case where $k=d$ and $f: X \to B$ is a family of 
Calabi-Yau manifolds.

\section{Cycle class maps from motivic cohomology into Deligne cohomology}
\label{chernclass}

Recall the following definition:

\begin{definition} A Calabi-Yau variety $X$ is a projective complex manifold of
dimension $d$ with the property that the canonical bundle 
$K_X=\Omega_X^d={\mathcal O}_X$ 
is trivial and $h^0(\Omega^i_X)=0$ for $1 \le i \le d-1$.
\end{definition} 

Why are Calabi-Yau manifolds interesting for our purpose?
The point is that the motive $h^d(X)$ is interesting, as shown first by
Griffiths \cite{G1} and later generalized by Voisin \cite{V1} in the case
where $d=3$. In fact, what we do in the following is equally interesting in the 
more general case where $K_X$ is trivial, e.g. in the case of abelian varieties. \\
There are cycle class maps defined by Bloch (see \cite{B3}) 
which can be realized by explicit Abel-Jacobi type integrals. There are formulas in
\cite{Go} for real Deligne cohomology, general formulas have recently been given 
in \cite{KLM} which we describe now. The maps
$$
c_{p,n}: CH^p(X,n) \longrightarrow H_{\mathcal D}^{2p-n}(X,{\mathbb Z}(p))
$$
were defined first by Bloch in \cite{B3} using Deligne cohomology with supports
for classical cycles and a spectral sequence construction. 
The abelian groups $ H_{\mathcal D}^{2p-n}(X,{\mathbb Z}(p))$ sit in exact
sequences
$$ 0 \to J^{p,n}(X) \to  H_{\mathcal D}^{2p-n}(X,{\mathbb Z}(p)) \to 
F^p H^{2p-n}(X,{\mathbb Z}(p)) \to 0,
$$ 
where 
$$
J^{p,n}=\frac{H^{2p-n-1}(X,{\mathbb C})}{F^p +H^{2p-n-1}(X,{\mathbb Z}(p))}
$$ 
are generalized intermediate Jacobians. Note that, if $n=0$ we get the
Griffiths intermediate Jacobians and for $p=1,n-1$ (resp.) the classical
Jacobian, resp. Albanese tori, which are Abelian varieties. Griffiths
intermediate Jacobians are not abelian varieties in general but they are 
non-degenerate compact tori. The generalized intermediate Jacobians are no 
longer compact, but they are still complex manifolds and vary holomorphically 
in families. If we restrict to cycles homologous to zero, we obtain Abel-Jacobi 
type maps \cite{B3,K,KLM} 
$$
c_{p,n}: CH^p(X,n)_{\rm hom} \longrightarrow J^{p,n}(X).
$$
For those we have explicit formulas found be M. Kerr \cite{K} 
(real versions can be found in \cite{Go}): consider an irreducible subvariety 
$W \subset X\times \square^{n}$, of codimension $p$, (using the cubical complex 
description), with coordinates $(z_{1},\ldots,z_{n})$, together with projections
$\pi_{1} : X\times \square^{n} \to X$ and $\pi_{2} : X\times \square^{n}
\to \square^{n}$. Let $\alpha \in F^{d-p+1}H^{2d-2p+n+1}(X,\C)$. 
One considers the current associated to $W$:

\begin{tiny}
$$ \alpha \mapsto
\frac{1}{(2\pi i)^{d-p+n}}\biggl[ \int_{W\backslash 
(W\cap \pi_{2}^{-1}([0,\infty]\times \square^{n-1})} \pi_{2}^{\ast}\big(
\log z_{1}d\log z_{2}\wedge\ldots\wedge d\log z_{n}\big)\wedge
\pi_{1}^{\ast}\alpha
$$
$$
-(2\pi i)\int_{\{W\cap \pi_{2}^{-1}([0,\infty]\times \square^{n-1})\}
\backslash \{W\cap \pi_{2}^{-1}([0,\infty]^{2}\times \square^{n-2})\}}
\pi_{2}^{\ast}\big(\log z_{2}d\log z_{3}\wedge\ldots\wedge 
d\log z_{n}\big)\wedge \pi_{1}^{\ast}\alpha
$$
$$
+\cdots + (-1)^{n-1}(2\pi i)^{n-1}\int_{W\cap \pi_{2}^{-1}
([0,\infty]^{n-1}\times \square^1)}\pi_{2}^{\ast}(\log z_{n})\wedge
\pi_{1}^{\ast}\alpha
+ (-1)^{n}(2\pi i)^{n}\int_{\Gamma}\pi^{\ast}_{1}\alpha  \biggr],
$$ 
\end{tiny}

where $\partial \Gamma = W\cap \pi_{2}^{-1} ([0,\infty]^{n})$ (the existence of 
$\Gamma$ follows from $W$ being homologous to zero, since the homology class
of $W$ is given by $W\cap \pi_{2}^{-1} ([0,\infty]^{n})$).
If $2p-n-1=d$ (i.e., $n=2p-d-1$), then we can truncate these maps and get
$$
\overline c_{p,n}: CH^p(X,n)_{\rm hom} \longrightarrow \frac{H^0(X,\Omega_X^d)^*}{
H_{d}(X,{\mathbb Z})}.
$$
The group on the right hand side is no longer a manifold since we
divide by a subgroup of rank bigger than the real dimension of
$H^0(X,\Omega_X^d)$, which is $2$ in our case. If we assume furthermore that
$p \le d$, then, given a cycle $W \in CH^p(X,2p-d-1)$, we have the formula
$$
\overline c_{p,n}(W)= (-1)^n \frac{1}{(2\pi i)^{d-p}} \cdot \int_\Gamma \omega_X,
$$
since the holomorphic $d$-form $\alpha$ vanishes on $W$, as
$\dim(W)=p-1 < d$. $\Gamma$ is as above and satisfies
$$  
\partial \Gamma= 
\sum_j (f^{(j)}_1,...,f^{(j)}_n)^{-1}([0,\infty]^n)=
\sum_j \bigcap_{\ell=1}^n (f^{(j)}_\ell)^{-1}([0,\infty]).
$$
Note that $\dim_{\mathbb R}(\Gamma)=2\dim(W)-n=d$ in this 
situation. If $p>d$ then the truncated Abel-Jacobi map 
involves more integrals. An interesting example is 
$CH^2_{\rm hom}(X,2)$ for $d=1$, i.e., where $X$ is an 
elliptic curve. More generally such examples are given by 
the groups $CH^{d+1}(X,d+1)$ for arbitrary $d$.

\section{Higher Chow groups of algebraic surfaces} 

Let $X$ be a complex projective surface and $Z \subseteq X$ a 
normal crossing divisor with complement $U=X \setminus |Z|$. Then 
we have an exact localization sequence
$$
\ldots \to CH^{r+1}(U,r+1) {\buildrel \partial \over \to} CH^r(Z,r) \to 
CH^{r+1}(X,r) \to \ldots 
$$
We are interested in the groups $CH^{r+1}(X,r)$ and, in view of this sequence, 
we are looking to have control over the group $CH^r(Z,r)$ which is of Milnor-type 
and hence combinatorially easy to describe. We also would like to give criteria
when classes in $CH^r(Z,r)$ come from $CH^{r+1}(U,r+1)$ via the coboundary map 
$\partial$. There is a spectral sequence
$$
E_1^{a,b}=CH^{a+n}(Z^{[-a]},1-b) \Longrightarrow CH^{n-1}(Z,-a-b)
$$
computing higher Chow groups of $Z$ from its smooth components $Z_i$.
Let $Z^{[1]}:=\coprod_i Z_i$ and $Z^{[2]}:=\coprod_{i<j} Z_i \cap Z_j$.
In particular for $r=1,2$, we have
$$
CH^r(Z^{[1]},r) \to CH^r(Z,r) \to CH^r(Z,r)^{\rm ind} \to 0,
$$
and 
$$
CH^r(Z,r)^{\rm ind}={\rm Ker}\left( CH^{r-1}(Z^{[2]},r-1) \to 
CH^{r}(Z^{[1]},r-1) \right).
$$
Hence 
$$
CH^1(Z,1)^{\rm ind}={\rm Ker}\left(\Z[Z^{[2]}] \to {\rm Pic}(Z^{[1]})\right), 
$$
and
$$
CH^2(Z,2)^{\rm ind}={\rm Ker}\left(CH^1(Z^{[2]},1) \to CH^2(Z^{[1]},1) \right). 
$$
For $r=1$, this observation had been the motivation to find complex projective 
surfaces with large motivic cohomology groups $CH^2(X,1)$ in \cite{SMS}:

\begin{theorem} \cite{SMS}  Let $X \subseteq \P^3$ be a very general quartic 
hypersurface containing a line. Then $CH^2(X,1)^{\rm ind}$ is not a torsion 
group, i.e., it contains elements
with no integer multiple which is decomposable.
\end{theorem} 

In \cite{Co}, Collino has independently found indecomposable cycles in $CH^2(X,1)$ 
on Jacobian surfaces. Related examples were found by Gordon and Lewis \cite{GL}, 
the authors \cite{dAMS}, Morihiko Saito \cite{S2}, Voisin and Oliva (both unpublished). 
In \cite{MSS} this method was further exploited in the case $r=2$ 
using some new techniques of Asakura and Saito:

\begin{theorem} \cite{SMS, MSS} Let $X \subseteq \P^3$ be a very general 
hypersurface of degree $d$ and let $Z=\cup Z_i$ be a NCD with $Z_i$ very 
general hypersurface sections of degrees $e_i$. Then \\
(1) If $d \ge 5$, $CH^2(U,2) \to CH^1(Z,1)^{\rm ind}$ is surjective. \\
(2) Assume $d \ge 6$ and that all triples $(e_i,e_j,e_k) \neq (1,1,2)$. Then
$CH^3(U,3) \to  CH^2(Z,2)^{\rm ind}$ is surjective as well. 
\end{theorem} 

Using the same methods one can also detect cycles
in $CH^2(X,1)$ and $CH^3(X,2)$. Here is a typical application:

\begin{example} Let $K$ be a general cubic polynomial in $\C[x_0,...,x_3]$. 
On the family 
$ X_{u,v}= \{F_{u,v}=x_0^5 + x_1x_2^4+x_2x_1^4+x_3^5 + ux_1^2x_2^3
+vx_0x_3 K=0 \}, \quad u,v \in \C$,
of quintic surfaces, there exist elements $Z_{u,v}$ in $CH^3(X_{u,v},2)$ 
such that, for $u,v \in \C$ very general, every integer multiple of 
these elements is not contained in the image of 
${\rm Pic}(X_{u,v}) \otimes K_2(\C)$. 
\end{example}

Note that the parameter $v$ is redundant since it is already contained in the 
coefficients of the cubic form $K$. However in the proof we fix 
$K=x_0^2x_1+x_0x_1x_2+x_0x_2^2+x_2^3+x_0x_1x_3$ and vary $u$ and 
$v$ in a 2-dimensional local parameter space to compute the infinitesimal data.
In an appendix to \cite{MSS}, Collino has constructed the following very 
interesting examples on K3 surfaces: 

\begin{theorem} \cite{MSS} On every (very) general quartic $K3$-surface $S$, 
there exists a 1--dimensional family  of elements $Z_t$ in $CH^3(S,2)$ such that, 
for $t$ very general, every integer multiple of these elements is not contained in
the image of ${\rm Pic}(S) \otimes K_2(\C)$. 
\end{theorem}

All these cycles originate from the existence
of smooth bielliptic hyperplane sections $C$ of genus $3$ on $S$, which means
that there exists a double cover $C \to E$ onto a smooth elliptic curve $E$.
Their existence is guaranteed by the following fact: 
on every general quartic surface $S$ there exists a
$1$-dimensional family of bielliptic curves $C_t$ such that the underlying
family $E_t$ of elliptic curves has varying $j$-invariant. \\
\ \\
We refer to section \ref{examples} for the construction of indecomposable elements
in $CH^2(X,1)$ of other K3-surfaces which can be detected using differential 
equations as in \cite{dAMS}.

\section{Poisson structures on generalized intermediate Jacobians}
\label{symplectic} 

Here we define the structure of a completely integrable Hamiltonian system
with in general non-compact leaves on the 
generalized intermediate Jacobians from section \ref{chernclass}. The normal 
functions associated to cycles in higher Chow groups are shown to be isotropic
subvarieties of the corresponding Poisson (i.e., degenerate symplectic) structure.
In the case of classical intermediate Jacobians this is a theorem of Donagi and 
Markman \cite{DM1}. See also Voisin's book \cite[page 12-14]{V2} and compare with
\cite{M} in the case of elliptic curves. 
Many famous classical examples of integrable systems can be found in 
\cite{knoerrer}.

\begin{definition}{\rm (cf. \cite{DM1})} A holomorphic symplectic structure on a 
complex manifold $M$ is given by a holomorphic (2,0)--form 
$\sigma \in H^0(X,\Omega^2_X)$ such that the contraction operator
$$
\rfloor \sigma: T_M \to \Omega^1_M
$$  
induces an isomorphism. 
\end{definition}
In particular, the complex dimension of $M$ is even. 
Under this isomorphism, the differential $df$ of every function $f$ corresponds to 
a vector field $X_f$. The flow of $X_f$ determines a dynamical system and one 
wants to know when it is integrable. 
\begin{definition} A complex Poisson manifold is a complex manifold with a bracket
$\{-,-\}$ on ${\mathcal C}^\infty$-functions which is a derivation in each 
argument and satisfies a Jacobi identity (Lie algebra structure).   
Every symplectic structure induces a Poisson structure on $M$ via
$\{f,g\}:=\sigma(X_f,X_g)$. On the other hand, any Poisson manifold determines a
global tensor $\Psi \in H^0(M,\Lambda^2 T_M)$ with the property 
$\Psi(df,dg)=\{f,g\}$, which in the symplectic case is just 
$(\rfloor \sigma)^{-1}$.  The rank of a Poisson manifold $M$ is given by the 
rank of $\Psi$. 
An analytically completely integrable Hamiltonian system is a set of pairwise 
Poisson commuting analytic functions $F_1,\ldots,F_m$ such that the analytic map 
$F=(F_1,\ldots,F_m): M \to B\subset \C^m$ is submersive with
rank equal to $\dim M-\frac{1}{2}{\rm rank}\Psi$ at every point. 
\end{definition}
Each of the functions $F_i$ is called a Hamiltonian function and the corresponding 
vector fields are called Hamiltonian vector fields. In this case
Liouville's theorem guarantees that the fibers of $F$ have a natural affine 
structure and the flow of the Hamiltonian fields $X_{F_i}$ is linear. Hence 
they are generalized tori and compact if and only if
$F$ is proper. A generalized torus here is a quotient of $\C^n$ by a discrete 
subgroup $\Z^g$ where $2g \le n$. We say that the Hamiltonian system is 
algebraically completely integrable, if $F: M \to B$ is a surjective algebraic 
morphism, where degenerate algebraic fibers are allowed. 
This property implies that such systems are integrable via generalized
Theta--functions. Here is one more important definition:
\begin{definition}
A subvariety $Z$ of a symplectic manifold $(M,\sigma)$ is 
Lagrangian, if the tangent space $T_{Z,z}$ to each generic point 
$z \in Z$ is isotropic and co-isotropic, i.e., maximally isotropic. 
If $(M,\Psi)$ is only Poisson, then we replace
the tangent bundle $T_M$ by the image $T$ of the induced map
$\rfloor \Psi: \Omega_M^1 \to T_M$ in the definition of Lagrangian subspace, 
see \cite[Def. 2.5.]{DM1}. 
\end{definition}

Note that the rank of $\Psi$ can vary, but there always exists a stratification 
of $M$ into Poisson strata of constant rank $r$ of $\Psi$, which are themselves 
foliated into symplectic leaves, i.e., on which there is an induced symplectic 
structure of rank $r$. \\
Let $f: X \to B$ a smooth and locally universal 
family of $d$-dimensional Calabi-Yau manifolds with smooth base $B$ 
and ${\mathcal H}^d$ be the flat vector bundle on $B$ associated 
to the $d$-th cohomology of the fibers. We also denote by 
${\mathcal H}_{\mathbb Z}^d \subset {\mathcal H}^d$ the local system of
integral cohomology classes. The set of Hodge bundles is given by holomorphic 
subbundles   
$$ 
{\mathcal F}^d \subseteq \ldots \subseteq {\mathcal F}^1 \subseteq 
{\mathcal F}^0={\mathcal H}^d.
$$
The bundle ${\mathcal F}^d$ is the same as the relative dualizing sheaf
$f_*\omega_{X/B}$ and is a line bundle since the fibers are Calabi-Yau manifolds. 
These vector bundles are equipped with the Hodge metric. In particular we
have metric isomorphisms for all $i$, 
$$ 
{\mathcal H}^d/{\mathcal F}^i \cong ({\mathcal F}^{d-i+1})^\vee,
$$
induced by the duality pairing 
$$
\langle-,-\rangle: {\mathcal H}^d \otimes {\mathcal H}^d \to  {\mathcal O}_B.
$$
The Gauss-Manin connection on ${\mathcal H}^d$ is given as 
$$ 
\nabla: {\mathcal H}^d \to {\mathcal H}^d \otimes \Omega_B^1. 
$$
Using base change, we now pass to $\tilde B$ which is the 
total space of the line bundle $f_*\omega_{X/B}$ with the zero-section
removed. All vector bundles on $B$ can be pulled back to $\tilde B$ along the 
natural projection map $p: \tilde B \to B$ and they are denoted by 
$\tilde {\mathcal H}^d$, $\tilde {\mathcal H}_{\mathbb Z}^d$ and 
$\tilde{\mathcal F}^i$ respectively. Assume
now that $2p-n-1=d$, not necessarily with $p \le d$. 
The family of generalized intermediate Jacobians ${\mathcal J}^{p,n} \to B$
pulls back to $\tilde {\mathcal J}^{p,n} \to \tilde B$. Note that
$\tilde{\mathcal F}^d$ has a tautological section $\eta$, which in each point of
$\tilde B$ picks up the holomorphic $d$-form defined by this point. The
Gau\ss-Manin connection pulled back to $\tilde B$ and applied to $\eta$
induces an isomorphism 
$$
\nabla_\bullet \eta: T_{\tilde B} {\buildrel \cong \over \longrightarrow} 
\tilde {\mathcal F}^{d-1}
$$ 
and hence, dually, an isomorphism between $\Omega^1_{\tilde B}$ and
$\tilde {\mathcal H}^d/\tilde {\mathcal F}^{2}$, see \cite[thm. 7.7]{DM1}. 
The proof of this fact is easy and consists of identifying the exact sequences 
$$
0 \to T_{\tilde B/B} \to T_{\tilde B} \to p^*T_B \to 0,
$$
and 
$$
0 \to \tilde{\mathcal F}^d \to \tilde{\mathcal F}^{d-1} 
\to \tilde{\mathcal F}^{d-1}/\tilde{\mathcal F}^d \to 0,
$$
via the map $\nabla_\bullet \eta$ and the identification of 
$T_{B,b}=H^1(X_b,T_{X_b})$ with $H^1(X_b,\Omega_{X_b}^{d-1})$, using the natural 
isomorphism $T_{X_b} \cong \Omega_{X_b}^{d-1}$. \\ 
Since the cotangent bundle on any complex manifold carries a natural
holomorphic symplectic structure, 
we conclude that $\tilde {\mathcal H}^d/\tilde {\mathcal F}^{2}$ 
carries a natural holomorphic symplectic structure, given by a 2-form
$d\theta$, where $\theta$ is the tautological 1-form on $\Omega^1_{\tilde B}$. 
We need the following description of $\theta$ from
\cite{V2}: at a point $(\tau,b)$ of 
$\tilde {\mathcal H}^d/\tilde {\mathcal F}^{2}$ over $b \in \tilde B$, we
have that 
$$
\theta_{(\tau,b)}=\pi^*\langle \tau, \nabla \eta \rangle,
$$ 
where $\pi: \tilde {\mathcal H}^d/\tilde {\mathcal F}^{2} \to \tilde B$ is the 
projection map. This fact follows immediately from the definition of the contact 
form on any cotangent bundle together with the isomorphism induced by 
$\nabla_\bullet \eta$.
The resulting $2$-form $d\theta$ is a closed and non-degenerate holomorphic 
$2$--form on $\tilde {\mathcal H}^d/\tilde {\mathcal F}^{2}$ and defines therefore 
a symplectic holomorphic structure such that the fibers of $\pi$
are easily seen to be Lagrangian as in \cite{V2,DM1}. \\ 
Now the universal covering space of $\tilde {\mathcal J}^{p,n}$
is given by $\tilde {\mathcal H}^d/\tilde {\mathcal F}^{p}$, which maps 
surjectively onto $\tilde {\mathcal H}^d/\tilde {\mathcal F}^{2}$ if $p \ge 2$. 
Therefore we can pull back $d\theta$ to 
$\tilde {\mathcal H}^d/\tilde {\mathcal F}^{p}$ to get a
Poisson structure. Note that the rank does not change under pullback, 
so that we cannot expect to get a symplectic structure anymore if $p>2$. 
However the fibers of the map to $\tilde B$ are still maximally isotropic
subspaces since they come from pullbacks of maximally isotropic subspaces. 

\begin{lem} $d\theta$ descends to a closed $2$-form $\sigma$ on 
$\tilde {\mathcal J}^{p,n}$. Furthermore all fibers of 
$\pi: \tilde {\mathcal J}^{p,n} \to \tilde B$ are maximally isotropic.  
\end{lem}
\begin{proof} As we remarked above, $\theta$ is 
given by the inner product $\langle -,\nabla \eta \rangle $
induced by the Hodge pairing on ${\mathcal H}^d$ and where $\eta$ is as above. Now, 
if two sections $\tau,\tau'$ of $\tilde {\mathcal H}^d/\tilde {\mathcal F}^p$ differ
by a flat section $\lambda \in \tilde {\mathcal H}_{\mathbb Z}^d$, then 
$\theta$ differs by $\pi^*\langle \lambda,\nabla \eta \rangle$. But since
$\lambda$ is flat, 
$$
d \langle \lambda,\nabla \eta \rangle = \langle \nabla \lambda,\nabla \eta
\rangle=0
$$
and $d\theta$ is therefore invariant. Let $\sigma$ be the induced form. The 
intermediate Jacobian fibers are Lagrangian by construction, since
they have this property for $p=2$.
\end{proof}

With slightly more work one can even show that $\sigma$ remains an exact form, 
see \cite{DM1}. For $p=2$, this structure will descend to an honest symplectic
structure on the total space $\tilde {\mathcal J}^{2,n} \to \tilde B$ if $d \le 3$ 
where $n=3-d$. In this case the fibers are generalized Lagrangian tori. 
The reader may wonder about the Hamiltonian functions: they only come up when we 
choose a level structure, i.e., a basis $(a_i)$ for the homology group 
$H_d(X_b,\Z)$, respecting the polarized structure, e.g.
a symplectic basis for $d$ odd. Then the fiber integrals
$$
F(b)=\int_{a_i} \eta(b)
$$
for $i=0,\ldots, b_d/2$ (assuming $d$ odd here) define pairwise commuting 
Hamiltonian functions, called the canonical coordinates in mathematical physics, 
see \cite{CdO1}. \\
On the other hand, if $p=1$, then by our assumption 
$2p-n-1=d$ we only get the possibility $d=1$ and $n=0$, which leads essentially 
to the consideration of $CH^1(X)$ for an elliptic curve $X$. \\
So far we have proved the first part of the theorem of
Donagi and Markman. Their second result says that (for $d=2p-1$ and $n=0$) 
all quasi-horizontal normal functions (e.g., arising from cycles 
in higher Chow groups which are defined over $B$ and hence over 
$\tilde B$) have Lagrangian subspaces as images under the 
Abel-Jacobi map. \\
A similar result holds for generalized intermediate Jacobians, but note that
if $p \ge 3$, then we only have a surjection 
$\tilde {\mathcal H}^d/\tilde {\mathcal F}^{p} \to 
\tilde {\mathcal H}^d/\tilde {\mathcal F}^{2}$   
and therefore we can only expect the image of the normal functions 
to be isotropic. We also refine the result to the case where the normal 
function is defined on an algebraic sublocus $B' \subset B$:
 
\begin{theorem} Assume that $p \le d=2p-n-1$ and $p \ge 2$. 
Let $\nu: B' \to  {\mathcal J}^{p,n}$ any quasi-horizontal normal function
defined on a (locally) closed subset of $B$, i.e., one that satisfies 
$\nabla \nu \in {\mathcal F}^{p-1} \otimes \Omega^1_B$. Then
the image of $\tilde \nu$ in $\tilde {\mathcal J}^{p,n}$ 
is an isotropic subvariety.
\end{theorem}
\begin{proof} The essential point is that $\tilde \nu$ is isotropic 
if and only if $\tilde \nu^* \theta$ is a closed form, because then
$d\theta$ will be zero along the image of $\tilde \nu$. To prove this, we
express again $\theta$ using $\langle -,\nabla \eta \rangle$ and we get
$$
d(\tilde \nu^* \theta)=d \langle \tilde \nu,\nabla \eta \rangle 
= d \langle \nabla \tilde \nu, \eta \rangle =0,
$$
since, by quasi-horizonality, we have 
$\nabla \tilde \nu \in \tilde {\mathcal F}^{p-1} \otimes \Omega^1_B$ and
the pairing $\langle \tilde {\mathcal F}^{p-1}, \eta \rangle $ already gives zero by
type reasons, since $\eta \in \tilde {\mathcal F}^{d}$ and $p-1 \ge 1$, i.e., $\eta$
pairs to zero with $\tilde {\mathcal F}^{p-1}$.  
\end{proof}
Note that a normal function may be multivalued. The theorem shows that
such normal functions defined on the whole of $B$ give rise to Lagrangian 
(maximally isotropic) subvarieties. We would also like to say that the projection 
of $\tilde \nu$ in $\tilde {\mathcal J}^{p,n}/\tilde {\mathcal F}^1$ is an 
isotropic subvariety, but this is not in general a manifold anymore.  
\begin{question} Do exotic symplectic structures like 
${\mathcal J}^{p,n}/{\mathcal F}^1$
have a rigorous mathematical meaning in another context? 
\end{question}

\section{More examples}
\label{examples}

As we remarked in section \ref{symplectic}, the most natural and 
interesting {\it time--independent} Hamiltonian systems occur in dimensions
$d \le 3$ when $p=2$. The case $d=1$ was treated in the paper of Manin \cite{M}: 
for the family of elliptic curves given by the equation $y^2=x(x-1)(x-t)$, one
has the inhomogeneous equation
$$
\left[ t(1-t)\frac{\partial^2}{\partial t^2}
+ (1-2t)\frac{\partial}{\partial t}-\frac{1}{4} \right]
\frac{dx}{y} = \frac{1}{2} \; d_{E/\P^1}\frac{y}{(x-t)^2},
$$
where the derivatives are understood in the covariant sense, i.e.,
$\frac{\partial x}{\partial t}=0$ since $x$ is a flat coordinate and 
$d_{E/\P^1} t=0$ by definition of the relative differential. In uniformized 
coordinates the non-linear equations take the equivalent form
$$
\frac{d^2}{d\tau^2} z(\tau)= \frac{1}{(2\pi i)^2} \sum_{2e=0} 
\alpha_e  {\mathfrak p}_z(z+e,\tau),
$$
where $e$ runs through all 2-torsion points and ${\mathfrak p}_z$ is the
derivative of the Weierstrass ${\mathfrak p}$-function. 
The Hamiltonian function governing the dynamical flow
of the {\it time-dependent} Hamiltonian system in \cite{M} is given by
$$
H=\frac{y^2}{2}-\frac{1}{(2\pi i)^2}\sum_{2e=0} \alpha_e {\mathfrak p}(z+e,\tau) 
$$
in $(y,z,\tau)$-space. It is very interesting to 
compare this to our approach and formulate a 
time-dependent version, see \cite{dAMS2}. \\
In dimension $3$ we have the work of Griffiths \cite{G4}. 
This case is only partly understood and we just
mention that there is a possibility that families of cycles in 
$CH^2_{\rm hom}(X)$ on Calabi-Yau 3-folds 
(e.g. differences of pairs of lines on quintic 3-folds) can be characterized by 
their inhomogenous Picard-Fuchs differential equations. 
The family of quintic 3-folds is discussed in 
\cite{CdO2,CdO3} in connection with arithmetic questions. \\
It remains to treat the case $d=2$: here we have started to 
give some examples related to Kummer surfaces
in \cite{dAMS}. Their equations were given by
$$
X_b=\{(x:y:z:w) \in \P^3 \mid xyz(x+y+z+bw)+w^4=0 \}.
$$  
This is a one-dimensional family of K3 surfaces (quartic) with generic 
Picard number $19$ (see \cite{NS} for this and other properties). 
Its mirror family is the family of all quartic K3 surfaces 
$$
f={\displaystyle \sum_{i=1}^4 x_i^4-t\prod_{i=1}^4 x_i}=0.
$$ 
let us work out the inhomogenous Picard-Fuchs equation in this example:
the algorithm of \cite{Mo1} produces the inhomogeneous equation
\begin{small}
$$
((t^4-256)(t\frac{\partial}{\partial t})^3w+2t^4(t\frac{\partial}{\partial t})^2
+\frac{7}{6}t^4(t\frac{\partial}{\partial t})-\frac{1}{6}t^4)\omega=
\frac{1}{3}dF[3]+\frac{1}{6}dF[2]+\frac{1}{6}dF[1],
$$
\end{small}
where
$$
\begin{array}{ll} 
F[3] =& \frac{1}{f^3} (  (-16x_1t^5x_2x_3^4x_4^4+64x_2^4t^4x_3^3x_4^3)dx_3dx_4 \\
  & +
(-x_1^3t^7x_2^2x_3^3x_4^2+64x_3^4t^4x_2^3x_4^3)dx_2dx_4 \\
  & +(-4x_1^2t^6x_2x_3^5x_4^2+64x_4^4t^4x_2^3x_3^3)dx_2dx_3 \\
  & +(-x_2^3t^7x_1^2x_3^3x_4^2+16x_3^5t^5x_2x_4^4)dx_1dx_4 \\
  & +(-4x_2^2t^6x_1x_3^5x_4^2+16x_4^5t^5x_2x_3^4)dx_1dx_3 \\
  & +(-4x_3^6t^6x_1x_2x_4^2+x_4^3t^7x_1^2x_2^2x_3^3)dx_1dx_2 ), \\
  & \\
  &
\end{array}
$$
$$
\begin{array}{ll}
F[2] =& \frac{1}{f^2}( -x_2(8t^5x_1x_3^4-2t^6x_1^2x_2x_3x_4)dx_3dx_4 \\ 
&+(x_1(-4t^5x_3x_4^4-8t^5x_1^4x_3)-x_3(8t^5x_1x_3^4-2t^6x_1^2x_2x_3x_4))dx_2dx_4 \\ 
 & +(-x_1^2t^6x_2x_3x_4^2-x_4(8t^5x_1x_3^4-2t^6x_1^2x_2x_3x_4))dx_2dx_3 \\ 
  & +x_2(-4t^5x_3x_4^4-8t^5x_1^4x_3)dx_1dx_4-x_2^2x_1t^6x_3x_4^2dx_1dx_3 \\ 
  & +(-x_3^2x_1t^6x_2x_4^2-x_4(-4t^5x_3x_4^4-8t^5x_1^4x_3))dx_1dx_2), \\
\end{array}
$$
and
$$
\begin{array}{ll}
F[1] =& \frac{1}{f}(  -2x_2t^5x_1dx_3dx_4-4x_1t^5x_3dx_2dx_4 \\
 &  -x_1t^5x_4dx_2dx_3  -2x_2t^5x_3dx_1dx_4 \\
 &  +x_2t^5x_4dx_1dx_3  +3x_3t^5x_4dx_1dx_2).
\end{array}
$$

\section{Outlook}
\label{outlook}

In this section we want to give an outlook for the reader and present some ideas 
which came up after the conference in Toronto. 
Assume again that $2p-n-1=d$ and $p \le d$. For every family $f:X \to B$ and 
every cycle ${\mathcal Z} \in CH^p(X,n)$ which is defined after restriction to 
every fiber, we get a complex function $\overline \nu$,
$$
\overline \nu(b)= \bar c_{p,n}({\mathcal Z}|X_b)=\int_{\Gamma_b} \omega_{X_b},
$$
induced by the truncated Abel-Jacobi map $\bar c_{p,n}$ from section 2 on each 
member of the family. The function $\overline \nu(b)$ is holomorphic on $B$ 
but has singularities along the boundary $\bar B \setminus B$
of a compactification $\bar B$. Using the notations from section 1, we have:
 
\begin{theorem} \label{extension} {\rm \cite{dAMS2}} For all symbols $\sigma$, 
$D_\sigma \overline \nu(b)$ has at most poles 
in the points in $\bar B \setminus B$, i.e., becomes a rational function 
on the compactification $\bar B$. 
\end{theorem}

\proof See \cite{dAMS2}. \endproof
The proof is a clarification and generalization of the result in 
\cite{dAMS} for the case of K3 surfaces. It uses Clemens' extension theory 
of normal functions together with the results of \cite{S1}. 
This result can be used to prove several consequences.
Assume that $f: X \to B$ is as above with $B$ and the cycle ${\mathcal Z}$
being defined over $\bar \Q$. Such a situation can for example be achieved by
spreading out a cycle on a generic fiber $X_\eta$ over the field obtained by the 
compositum of its field of definition and the function field of $\eta$. In other 
words all transcendental elements in the equations of $X_\eta$ and 
${\mathcal Z}$ occur in the coordinates of $B$.
Then there is a unique choice of relative canonical
$d$-form by spreading out also $\omega$ from $X_\eta$.

\begin{theorem} In this situation, the rational function $g$ has coefficients in
$\bar \Q$.
\end{theorem}

\proof See \cite{dAMS2}.
\endproof
Such questions can be used to reprove the countability result
\cite[Cor. 3.3]{SMS}. \\
\ \\
In the paper of Manin \cite{M}, the general formalism
of {\it $\mu$--equations} in the context of families of algebraic curves was
established. This was used by him already in the solution of the functional Mordell
conjecture. Let us explain how to generalize this and assume that
we can solve $D_\sigma \omega_X=d \beta$ for some
$(d-1)$-form $\beta$, where $D_\sigma$ is one of the generators of the
differential ideal of Picard-Fuchs equations, see section \ref{picardfuchs}.
If we have a fibration $\Gamma \to B$ where $\Gamma$ is a relative singular chain
of real dimension $d$, then we are looking at the equation
$$
D_\sigma \int_\Gamma \omega = \int_\gamma i^* \beta + T,
$$
for some relative $(d-1)$-form $\beta$ and where $i: \gamma \hookrightarrow X$
is the inclusion of $\gamma=\partial \Gamma$ in the total space $X$.
The term $T$ arises from differentiating
the boundary of $\Gamma$. This formalism is a
direct generalization of \cite[formula (1.5)]{M}.
The connection between non-linear equations and inhomogenous
Picard-Fuchs equations is a direct consequence of these
$\mu$-equations.
In \cite{dAMS2} we give more examples of non-linear equations
arising from higher dimensional examples. \\
\ \\
{\it Acknowledgements:} We would like to thank the organizers of the workshop
on ``Arithmetic, Geometry and Physics around Calabi-Yau Varieties and Mirror
Symmetry'' and the staff of the Fields Institute for the opportunity to
present our research and visit the Fields Institute.
We are grateful for the patience of the editors and would like to thank
Xavier Gomez-Mont, Mark Green, Phillip Griffiths, James Lewis, Matthew Kerr and
Noriko Yui for several discussions. We thank
Jan Nagel and Morihiko Saito for help on normal functions.
This paper is dedicated to Andrei Tyurin who has encouraged us to work in this
direction and explained to us his ideas about mirror symmetry.

\bibliographystyle{plain}
\renewcommand\refname{References}

\end{document}